\documentclass[11pt]{amsart}
\usepackage{amsmath, amsthm, amssymb}
\DeclareMathOperator{\Span}{span}
\DeclareMathOperator{\dist}{dist}

\newtheorem{thm}{Theorem}[section]
\newtheorem{lemma}[thm]{Lemma}
\newtheorem{cor}[thm]{Corollary}
\newtheorem{prop}[thm]{Proposition}
\newtheorem{example}[thm]{Example}
\newtheorem{question}[thm]{Question}
\newtheorem{define}[thm]{Definition}

\theoremstyle{remark}
\newtheorem*{remark}{Remark}

\author{Steven F. Bellenot}
\title{Skipped Blocking and other Decompositions in Banach spaces}

\begin{document}

\begin{abstract}
Necessary and sufficient conditions are given for when a sequence
of finite dimensional subspaces $(X_n)$ can be blocked to be
a skipped blocking decompositon ({\sf SBD}).
These are very similar
to known results about blocking of biorthogonal sequences.
A separable space
$X$ has {\sf PCP}, if and only if, every norming decomposition 
$(X_n)$ can be blocked
to be a boundedly complete {\sf SBD}. Every
boundedly complete {\sf SBD} is a {\sf JT}-decomposition.
\end{abstract}

\maketitle
\markboth{Steven F. Bellenot}{Skipped Blocking and other Decompositions}

\section{Introduction}

Skipped blocking decompositions ({\sf SBD}), collections of finite
dimension subspaces $(X_n)$ of a Banach space $X$ with additional properties
(Definition \ref{sbddn}), were introduced by Bourgain
and Rosenthal \cite{BR} to explore {\sf RNP}, the Radon-Nikodym Property,
and the weaker {\sf PCP}, the point of continuity property.
The standard Mazur product construction (Proposition \ref{Mazur}),
basically a global gliding hump,
shows that every separable Banach space $X$ has an {\sf SBD}.
This construction has a great deal of flexibility. Properties
of the underlining space $X$ can often be
transfered to the constructed sequence.
For example, Bourgain and Rosenthal \cite{BR} showed the existence of a 
boundedly complete {\sf SBD} (Definition \ref{bcsbddn}) implied {\sf PCP}.
Ghoussoub and Maurey \cite{GM1} using results form Edgar and Wheeler
\cite{EW} showed a converse, separable spaces with {\sf PCP} have a 
boundedly complete {\sf SBD}. Other examples are in \cite{R}, \cite{F},
\cite{B1} and \cite{B2}.
Having a boundedly complete {\sf SBD} implies the existence
of decompositions with stronger properties. Ghoussoub, Maurey and
Schachermayer \cite{GMS}
showed such spaces must
have a {\sf JT}-decompostion (Definition \ref{JTdn}).
Like many techniques, {\sf SBD} existed in the literature before the
technique was named. For example, the proof that each separable $X$ has
a subspace $Y$ so that both $Y$ and $X/Y$ have {\sf FDD} \cite{JR}
(see \cite{LT} page 48), is essentially blocking a norming biorthogonal
sequence into a {\sf SBD}.

Our unifying theme is the question: ``can one strengthen the {\em existence} of
a (nice) decomposition to the {\em universal} all decompositions
(are nice)?''
There is a trivial obstruction,
one might need to block the decomposition, the
sequence $(X_n)$, in
order to make it skipped blocking or make the 
the decompositon boundedly complete.
We introduce the notion of {\sf DDD} (Definition \ref{sbddn})
which are the sequences $(X_n)$ that might be blockable to be a {\sf SBD}. The
{\sf DDD} property is studied without a name in \cite{BR} and called
a {\em decomposition} with no adjectives in \cite{R}. Also {\sf DDD}
generalize the well known notation of a biorthogonal sequence. In fact,
{\sf DDD} are nothing more than blockings of biorthogonal sequences.
Our general
question becomes: ``when can every {\sf DDD} (with perhaps additional
properties) be blocked to be a {\sf SBD} (perhaps with additional properties)?''

The ordering of the sequence $(X_n)$ doesn't matter in the
definition of {\sf DDD}. Thus the collection of {\sf DDD} include
permutations of conditional basic sequences, which can be very ill behaved.
Theorem \ref{sbdiffnorming} states that
{\sf DDD} can be blocked to be a {\sf SBD} exactly when the 
predecomposition space $Y \subset X^*$ is $c$-norming for some $c$.
One can renorm $X$ to improve to make the $c$-norming constant $1$
(Proprosition \ref{renorm}). 
The space $Y$ is independent of the ordering on $(X_n)$.
If $(X_n)$ is {\sf DDD}, then the dual
system $(X_n^*)$ can be blocked to be a {\sf SBD} (Proposition \ref{dualddd}).

A number of examples are given that illustrate these results.
Example \ref{badddd} shows not every {\sf DDD} can be blocked to
be a {\sf SBD}.
Example \ref{trigbasis} shows that the dual sequence $(X^*_n)$ of a
{\sf DDD} need not be a {\sf SBD}.
Example \ref{J} is a permutation of
a basis for James space $J$ and provides a {\sf SBD} that has badly 
behaved partial sum projections no matter how the decomposition is blocked.
Example \ref{J2} shows the subspace $X_n$ can be far from the 
the quotient $X/[X_m]_{m\not=n}$.

We have two complete solutions for the boundedly complete {\sf SBD} case,
or equivalently for separable spaces $X$ with {\sf PCP}.
The point of continuity property, {\sf PCP},
states every bounded set has a point of weak to norm continuity.
Theorem \ref{allbc} shows every norming {\sf DDD} 
in a space with {\sf PCP} can be blocked to be
a boundedly complete {\sf SBD}.

A {\sf JT}-decomposition (Definition \ref{JTdn}) is a boundedly complete
skipped
decomposition with additional properties like those of the predual
of {\sf JT}, James tree space.
The fact that {\sf PCP} implies the existence of a {\sf JT}-decomposition
was proved in \cite{GMS} using results from sequence of previous papers
\cite{GM1}, \cite{GM2}. Their construction of {\sf JT}-decomposition 
added more conditions to the Mazur product construction from the
earlier boundedly complete {\sf SBD} construction.
We strengthen this result by showing that each boundedly
complete {\sf SBD} is already
a $c$-norming {\sf JT}-decomposition (Theorem \ref{bcdJT}) for some
$c < \infty$.
Theorem \ref{climax} shows
if $X$ has the {\sf PCP} then every {\sf DDD} with
a $c$-norming predecompostion spae can be blocked to be
a $c$-norming {\sf JT}-decomposition and the space $X$ can renormed
so that the blocking is a $1$-norming {\sf JT}-decomposition.

The name {\sf SBD} and its adjectives are somewhat unwieldy. The
name boundedly complete {\sf SBD} applies property ``boundedly complete''
only to skipped subsequences and not to the global decomposition.
The term {\sf DDD} doesn't demystify a defining-phrase or
stand for anything, but we wanted to reserve ``decomposition'' with
no adjective for informal use.

The author would like to acknowledge the help of a
referee of an earler version of the paper. For both some
connections with biorthogonal sequences and for the examples
about total vs norming vs $c$-norming. These appear at the end
of Section 3, starting with Proposition \ref{anyy}.

\section{Notation, Preliminaries, {\sf DDD} and {\sf SBD} }

We start the notation about {\sf DDD} and
skipped-blocking decompositions, {\sf SBD} in this section. Besides 
notation, we explored the theory with simple examples and
observations some of which are known. 
Also in this 
section is a well-known preliminary proposition. Every paper on 
skipped-blocking decompositions seems to have a proof based on the 
Mazur product construction and Proposition \ref{Mazur} is ours.

Our notation generally follows \cite{LT} or \cite{JL}, the first chapter
of \cite{JLE}.  In particular, 
$[X_n]^k_{n=1}$ is the closed linear span of $\cup^k_{n=1}X_n$ and 
$[X_n]=[X_n]^{\infty}_{n=1}$.   If $m\leq k$ are integers, we will 
write $X[m,k]$ for $[X_n]^k_{n=m}$ and $X[m,\infty)$ for $[X_n]_{n=m}^\infty$.
If $(m(i))$ is a strictly increasing integer sequence with $m(0) = 0$,
then we will say $(X[m(i-1)+1,m(i)])$ is a blocking of $(X_n)$. We
use the dual pair notation $\left\langle x, y \right\rangle$ for
$y(x)$ or $x(y)$.

\begin{define}
\label{sbddn}
We will say that $(X_n)$, a sequence of finite dimensional 
subspaces, is a {\em skipped-blocking decomposition} ({\sf SBD}) for a 
Banach space $X$ provided (1)--(3) hold. If only (1) and (2) hold
we will say $(X_n)$ is a {\sf DDD}.
\end{define}

\begin{enumerate}
\item[(1)]$X=[X_n]$. Sometimes this property is called {\em total}.

\item[(2)]For each $n$, $X_n\cap[X_m]_{m\not= n}=\{0\}$. Sometimes
this property is called {\em minimal}.

\item[(3)]For sequences $(n(i))$ and $(m(i))$ with 
$n(i)<m(i)+1<n(i+1)$ $(X[n(i),m(i)])^{\infty}_{i=1}$, is {\sf FDD} for 
$[X[n(i),m(i)]]$.
\end{enumerate}

Given an {\sf DDD} $(X_n)$ we define the projections $p_n$, $P_n$ and
$R_n$ as follows:

\begin{enumerate}
\item[(4)]The projection $p_n:X\rightarrow X$ with kernel $[X_m]_{m\not= 
n}$ and range $X_n$.

\item[(5)]The projection $P_n:X\rightarrow X$ given by 
$P_n=\sum\nolimits^n_{i=1}p_i$.

\item[(6)]The projection $R_n:[X_m]_{m\not= n}\rightarrow[X_m]_{m\not= n}$ 
which is the restriction of $P_{n-1}$ or $P_n$. The more general
skipped projections $R_{m,k}$ for $m \le k$, which is the restriction of
$P_m$ on the space $X[1,m-1]\oplus X[k+1,\infty)$. We have
$\Vert R_{n,j} \Vert \le \Vert R_{m,k} \Vert$ whenever $n\le m \le k \le j$.

\item[(7)]The constants $K=\sup\Vert R_n\Vert$,
$K_{\infty}=\lim\sup\Vert R_n\Vert$,
and $K_{\infty,\infty} = \lim_m \lim_k \Vert R_{m,k} \Vert$
Equation (3) is equivalent to $K<\infty$.  
Equation (2) implies that the projections in (4)--(6) are 
bounded. The monotone estimate in (6) implies that the limit
$K_{\infty,\infty}$ exists.
\end{enumerate}

\begin{define}
We will call the constant $K$ in (7), the {\em 
{\sf SBD}-constant}, and the constant $K_{\infty}$, the {\em asymptotic 
{\sf SBD}-constant}. Sometimes the constants of a {\sf SBD} can
be improved by blocking, the constant $K_{\infty,\infty}$ is the
{\em limiting asymptotic constant}, a {\sf DDD} can be blocked to be a
{\sf SBD}, if and only if, $K_{\infty,\infty} < \infty$. 
\end{define}

\begin{remark}[]
\label{notfdd}
The sequence $(X_n)$ is an {\sf FDD} exactly if the projections $(P_n)$
are uniformly bounded which would imply that the projections $(p_n)$ are
also uniformly bounded. Conversely, since $P_n = R_{n+1}(I-p_{n+1})$,
if $(p_n)$ are uniformly bounded and $(X_n)$ is a {\sf SBD},
then $(X_n)$ is an {\sf FDD}.
On the other hand, the principal of uniform
boundedness says if $\Vert P_n \Vert$ is unbounded, then there is
an $x \in X$ with $\Vert P_n x \Vert$ unbounded. We will see later
(Example \ref{J}), there are {\sf SBD} $(X_n)$ and $x$ where no
subsequence of $(P_n x)$ is bounded. The projections $(P_n x)$ can be very
far from $x$.
\end{remark}

\begin{example}
\label{badddd}
A {\sf DDD} of one-dimensional subspaces in Hilbert space that is not
a {\sf SBD}, nor can it be blocked to be a {\sf SBD}.
\end{example}

\begin{proof}[Construction]
The sequence $(X_n)$ where $X_n = [e_1 + e_{n+1}/n]$ in Hilbert space with
orthonormal basis $(e_n)$ satisfies both (1) and (2) but not (3),
so it is a {\sf DDD} but not a {\sf SBD}. Since $e_1+ e_{n+1}/n \rightarrow
e_1$,
the sets $[X_m]_{m\not=n} = [e_m]_{m\not= n+1}$ and
so the space $X^*_n = [e_{n+1}]$ (see (9) below).
No blocking of $(X_n)$ is a {\sf SBD}, but $(X^*_n)$
is a {\sf FDD} for its closed linear span. Since
$\cap_n X[n,\infty) = [e_1] \not= \{0\}$, the span of $(X^*_n)$
is not all of the dual.
Thus the {\sf DDD} is not separating (Definition  \ref{normdn}).
The projections $p_n(x) = e^*_{n+1}(x)(ne_1+e_{n+1})$ and hence
$p_n(e_1) = 0$ for all $n$. On the other hand $x=\sum e_{n+1}/n$
has $p_n(x) = e_1+e_{n+1}/n$ and $\Vert P_n x \Vert > n$.
Having $\cap_n X[n,\infty) = [e_1] \not= \{0\}$
is the only way a {\sf DDD} can fail to be a {\sf SBD} in Hilbert
space.
\end{proof}

\begin{example}
A {\sf SBD} of one-dimensional subspaces in Hilbert space that is not
a {\sf FDD}.
\end{example}

\begin{proof}[Construction]
Let $(e_n)$ be an orthonormal basis, let $X_n$ be
the one dimensional $[e_n]$ when $n$ is odd and the one dimensional
$[e_{n-1}+e_n/n]$ when $n$ is even. Eventually, $\Vert p_{2n} \Vert = 2n$,
and $\Vert R_n\Vert = 1$.  Thus $(X_n)$ is not a {\sf FDD} but is
a {\sf SBD} with constant one. This also shows even when the $X_n$ are
one dimensional in a {\sf SBD} there need not be a bound on $\Vert p_n\Vert$.
This is, of course, the standard example that can be found in many places.
Clearly the blocking given by $m(i) = 2i$ improves this {\sf SBD} to
a {\sf FDD}.
\end{proof}

\begin{prop}
\label{Mazur}
If $X$ is separable, then $X$ has a {\sf SBD} $(X_n)$ 
whose asymptotic constant is one.
\end{prop}

\begin{proof}[Outline of Proof]
Let $\varepsilon_n>0$ so that $\prod(1+\varepsilon_n)<\infty$.  Let $(x_n)$ be 
dense in $X$.
Let $X_1=[x_1]$ and let $W_0$ be a finite set of norm one 
elements of $X^*$ so that $X=X_1\oplus W^{\perp}_0$
Inductively pick 
$(X_n)\subset X$ and finite subsets $W_n$ of the unit sphere of 
$X^*$ so that

\begin{enumerate}
\item[(i)]$W_n\supset W_{n-1}$.

\item[(ii)]For $x\in X[1,n]$
$\Vert x\Vert\leq(1+\varepsilon_n)\sup\{|\langle x,y\rangle| : y \in W_n\}$.

\item[(iii)]If $x_{n+1}=u+v$ with $u\in X[1,n]$
$v\in W^{\perp}_{n-1}$ and $v\not= 0$, then there is $y\in W_n$ 
with $\langle v, y\rangle\not=0$.

\item[(iv)]$X_{n+1}\subset W^{\perp}_{n-1}$ so that both
$W^{\perp}_{n- 1}=X_{n+1}\oplus W^{\perp}_n$ and $v\in X_{n+1}$.  We have 
$X=X[1,n+1]\oplus W^{\perp}_n$.
\end{enumerate}

The proof is now clear.  But to belabor the point, note (ii) 
implies $\Vert R_{n+1}\Vert\leq 1+\varepsilon_n$ and (iii) implies $X=[X_n]$.
\end{proof}

\begin{example}
A {\sf SBD} which cannot be blocked to be a {\sf FDD}.
\end{example}

\begin{proof}[Construction]
Let $X$ be a Banach space without the Approximation Property, so in particular
it cannot have a {\sf FDD}. Since $X$ has a {\sf SBD} by
Proposition \ref{Mazur}, this {\sf SBD} of $X$ cannot be blocked to be
a {\sf FDD}.
\end{proof}

\begin{prop}
\label{condbasis}
For any basis $(e_i)$ of $X$ and and permutation $\pi$ the
sequence of one dimensional spaces defined by $X_n = [e_{\pi(j)}]_{j=n}^n$
is a {\sf DDD}. Furthermore,
$(X_n)$ is {\sf SBD}, if and only if, $(e_{\pi(n)})$ is a basis.
\end{prop}

\begin{proof}
The definition of {\sf DDD} is invariant under permutations. Since
$(e_i)$ is a basis, $x_n = e^*_{\pi(n)}(x)e_{\pi(n)} = p_n(x)$ is uniformly
bounded in norm.
The remark before Example \ref{badddd}
shows if $(X_n)$ is a {\sf SBD} then $(X_n)$
is an {\sf FDD}, and hence $(e_{\pi(n)})$ is a basis. The converse is formal.
\end{proof}

\section{Duality and the Predecomposition Space}

Continuing the list of notation, given an {\sf DDD} $(X_n)$ we define
the quotients $Z_n$ and $Z_{m,k}$ with quotient maps $q_n$, and
$q_{m,k}$ and subspaces $X_n^*$ and $Y$ (the predecomposition space) in
the dual, $X^*$. The quotient maps yield good estimates, unlike the
projects $(P_n)$ or $(p_n)$ as Example \ref{J} and \ref{J2} show. 

\begin{enumerate}
\item[(8)]The quotient map $q_n:X\rightarrow Z_n=X/[X_m]_{m\not= n}$.
And the more general quotients $q_{m,k}:X\rightarrow
Z_{m,k} = X/(X[1,m-1]\oplus X[k+1,\infty))$.
We
have $\Vert q_{m,k}(x) \Vert \ge \Vert q_{n,j}(x) \Vert$ whenever
$m \le n \le j \le k$.

\item[(9)]The subspaces $X^*_n=[X_m]^{\perp}_{m\not= n}\subset X^*$.
Clearly $X^*_n$ is isometric to $Z_n^*$ via the injection $q_n^*$. Also
$X^*_n$ is the range of the projection $p_n^*$. Finally,
$Z^*_{m,k} = [X^*_n]_{m\le n\le k}$.

\item[(10)]$X^{*\perp}_n=[X_m]_{m\not= n}$, for each $n$.

\item[(11)]$X^*_n\cap[X^*_m]_{m\not= n}=\{0\}$, for each $n$.

\item[(12)]$p^*_n$ is a projection on $X^*$ (and $Y$) with range $X^*_n$ and 
kernel $X^{\perp}_n\supset[X^*_m]_{m\not= n}$.
\end{enumerate}

\begin{remark}[]
Each {\sf DDD} $(X_n)$ is a blocking of a biorthogonal sequence $(x_i, x^*_i)$
obtained by combining  the biorthogonal sequences for the finite
dimensional space pairs $(X_n, X^*_n)$. The converse is also clear, each
grouping of a biorthogonal sequence yields a {\sf DDD}.
\end{remark}

\begin{define}
Given a {\sf DDD} $(X_n)$ let $Y_n = X^*_n$ and let
$Y = [Y_n] = [X_n^*] \subset X^*$.
We will say $Y$ is the predecomposition space of $(X_n)$.
Clearly $(Y_n)$ is a {\sf DDD} for $Y$ by equation (11). Note
that (9) implies that any blocking of $(X_n)$ has the same
predecomposition space $Y$.
\end{define}

\begin{prop}
\label{dualddd}
Any {\sf DDD} $(X_n)$ can be blocked to $(X[m(i-1)+1,m(i)])$
so that resulting dual {\sf DDD} $(Y[m(i-1)+1,m(i)])$
is a {\sf SBD} for the predecomposition
space $Y$.
\end{prop}

\begin{proof}
Suppose $(X_n)$ is a {\sf DDD}.
Let $\varepsilon > 0$ be given with $\varepsilon < 1$.
Once $m(k)$ is selected, consider the quotient $Z = Z_{1,m(k)}$
and subspace of $Y$ given by $W = q^*_{1,m(k)}(Z^*)$. We can find
a finite subsets $\{z_i\}_1^N$ and $\{w_i\}_1^N$ of the unit sphere
of $Z$ and $W$ so that $\left\langle z_i, w_i\right\rangle = 1$
and $\{w_i\}_1^N$ is an
$\varepsilon$-net for the sphere of $W$. Since $\Span \cup X_n$ is dense
we can find $x_i$ so that $1 \le \Vert x_i \Vert < 1+\varepsilon$ and
$q_{1,m(k)}(x_i) = z_i$. Finally select $m(k+1)$ so that 
$\{x_i\}_1^N \subset X[1,m(k+1)]$.

So if $y \in [X^*_n]_{n>m(k+1)}$ and $w$ is in  the sphere of $W$,
then by the construction above we can find $w_i$ with
$\Vert w-w_i\Vert < \varepsilon$ and clearly,
$y(x_i) = 0$. Hence $\Vert w+y \Vert \ge \Vert w_i+y\Vert - \Vert w-w_i\Vert$
and $\Vert w_i + y \Vert
\ge |(w_i+y)(x_i)|/\Vert x_i\Vert \ge 1/(1+\varepsilon)$. So we
have $\Vert w + y \Vert \ge (1-\varepsilon)(1+\varepsilon)^{-1} \Vert w\Vert$
and the $R_{m(k)+1,m(k+1)}$
projection onto $[X^*_n]_{n\le m(k)}$ with kernel $[X^*_n]_{n>m(k+1)}$
has norm bounded by $(1+\varepsilon)/(1-\varepsilon)$.
\end{proof}

\begin{question} 
\label{q3}
If $(X_n)$ is a {\sf SBD}, is the dual {\sf DDD} $(X^*_n)$ already a {\sf SBD}?
\end{question}

\begin{example}
\label{trigbasis}
A {\sf DDD} $(X_n)$ for $\ell_2 $, where $(X^*_n)$ is not a {\sf SBD} for
$Y = [X^*_n] = \ell_2$.
\end{example}

\begin{proof}[Construction]
Let $(e_i)$ be a conditional basis for $\ell_2$ and let $\pi$ be a permutation
so that $(e_{\pi(n)})$ is not a basis. Thus the coefficient functionals
$(e^*_{\pi(n)})$ are also not a basis. Letting $X_n = [e_{\pi(j)}]_{j=n}^n$,
Proposition
\ref{condbasis} says $(X^*_n)$ is not a {\sf SBD}.
\end{proof}

\begin{cor}
\label{condsbd}
For any basis $(e_n)$ and permutation $\pi$, the sequence $(e_{\pi(n)})$
can be blocked to be a {\sf SBD}.
\end{cor}

\begin{proof}
Apply Proposition \ref{dualddd} to the coefficient functions $(e^*_{\pi(n)})$.
\end{proof}

\begin{lemma}
\label{estimates}
Suppose $(X_n)$ is a {\sf DDD}, then the following estimates hold:
\begin{enumerate}
\item[(a)] If $w \in X[1,k-1]$ and $j \ge k$,
then $\Vert w \Vert \le \Vert R_{k,j}\Vert \Vert q_{1,j}(w)\Vert$ and
in particular,
$\Vert w \Vert \le \Vert R_k\Vert \Vert q_{1,k}(w)\Vert$.

\item[(b)] If $x \in X$, and $\delta_k = \dist(x, X[1,k-1])$ then
$\Vert x \Vert \le \Vert R_k \Vert\Vert P_k x\Vert + \delta_k(1 + \Vert R_k\Vert)$.

\item[(c)] If $(X_n)$ is a {\sf SBD}, then for all $x \in X$, 
$$\Vert x \Vert \le K_\infty\limsup \Vert P_n x \Vert \;\mathrm{and}\;
\Vert x \Vert \le K \liminf \Vert P_n x \Vert$$.

\item[(d)] If $w \in X[1,k-1]$, and $j \ge k$,
then there is $y_j \in [X^*_i]_1^j$ with
$\Vert y_j \Vert = 1$ and
$\Vert w \Vert \le \Vert R_{k,j}\Vert |\langle w, y_j \rangle |$.

\item[(e)] The predecomposition space $Y$ separates points in $X$, if
and only if, $\cap_n X[n,\infty) = \{0\}$.

\end{enumerate}
\end{lemma}

\begin{proof}
Parts of this proof are essentially stolen from the proofs of
Lemma I.15 and Lemma I.12 of \cite{GMS}.

Let $w \in X[1,k-1]$ and  
let $z \in X[j+1,\infty)$, and note $\Vert w \Vert = \Vert R_{k,j}(w+z) \Vert
\le \Vert R_{k,j}\Vert \Vert  w+z\Vert$.
Hence $\Vert w \Vert \le \Vert R_{k,j}\Vert  \Vert q_{1,j}(w) \Vert$
which proves (a).

Let $x \in X$, $\varepsilon > 0$ and $\delta_k = \dist(x, X[1,k-1])$.
Find $w \in X[1,k-1]$ with $\Vert w - x \Vert < \delta_k + \varepsilon$.
We have $w = P_k w$, 
$q_{1,k}(x) = q_{1,k}(P_k x)$ and
hence
$\Vert q_{1,k}(w) - q_{1,k}(P_k x)\Vert = \Vert q_{1,k}(w-x)\Vert
\le  \delta_k+\varepsilon$.
Putting the pieces together, 
\begin{align*}
\Vert q_{1,k}(P_k x) \Vert
+ \delta_k+\varepsilon
&\ge \Vert q_{1,k}(w) \Vert \ge \Vert w \Vert/\Vert R_k\Vert
\\
\Vert R_k\Vert \Vert P_k x \Vert
+ \Vert R_k\Vert (\delta_k + \varepsilon)
&\ge \Vert w \Vert \ge \Vert x \Vert-\delta_k -\varepsilon
\end{align*}
Which completes (b).

If $(X_n)$ is not a {\sf SBD}, then the inequalities in (c) are 
ambiguous, as Example \ref{badddd} shows that $K = K_\infty = \infty$
and $P_k x = 0$ for some $x$ and all $k$. The usual convention
$0\cdot \infty = 0$ gives the wrong result for general {\sf DDD}.
Since $\delta_k$ is monotonically decreasing to zero we can derive the
first part of (c) from
(b) by picking $k$ with $\Vert R_k \Vert \approx \limsup \Vert R_k\Vert$.
We derive the second part of (c) from (b) by picking $k$ with
$\Vert P_k x \Vert \approx \liminf \Vert P_k x \Vert$.

For (d) one uses (a)
and the fact $Z^*_{1,j} = [X^*_i]_1^j$ isometrically. So there is a
$y_j \in Z^*_{1,j}$ with $\Vert  y_j \Vert = 1$ and
$\langle w, y_j \rangle = \Vert q_{1,j}(w) \Vert$.

The statement (e) is almost immediate. If $x \not\in X[n+1,\infty)$, then
$q_{1,n}(x) \not= 0$. So there is a $x^* \in Z^*_{1,n} \subset Y$
with $\langle x, x^* \rangle \not= 0$. Thus $Y$ separates the points of
$X$ when $\cap_n X[n,\infty) = \{0\}$. Conversely, if there is a non-zero
 $x \in \cap_n X[n,\infty)$
then $\langle x, x^* \rangle = 0$ for $x^* \in X^*_m$ and any $m$. So
$x$ is zero on a dense subset of $Y$. Thus $Y$ doesn't separate the point
$x$ from $0$.
\end{proof}

\begin{define}
Let $c < \infty$,
space $Y \subset X^*$ is said to {$c$-norms} $X$, if all $x \in X$,
$$c^{-1}\Vert x \Vert  \le
\Vert x \Vert_Y =\sup\{|\langle x, y \rangle|: y \in Y, \Vert y \Vert\le 1\}
\le \Vert x \Vert$$
\end{define}

\begin{define}
\label{normdn}
A {\sf DDD} is said to be {\em separating}, if the predecompostion space
$Y$ separates the points of $X$ or equivalently if $\cap X[n,\infty) = \{0\}$.

A {\sf DDD} is said to be norming, if the predecomposition space $Y$ is
$c$-norming for some $c < \infty$.
\end{define}

\begin{thm}
\label{sbdiffnorming}
If $(X_n)$ is a {\sf DDD} then the predecomposition space $Y$ $c$-norms $X$
for $c = K_{\infty,\infty}$. Conversely if for some $c<\infty$, $Y$ $c$-norms
$X$ then the {\sf DDD} can be blocked to be a {\sf SBD} with asymptotic
constant $\le c$.
\end{thm}
\begin{proof}
Let $\varepsilon > 0$ be given.
If $K_{\infty,\infty} = \infty$ there is nothing to prove. Otherwise,
by blocking $(X_n)$ by $m(i)$ we can assume
$\Vert R_{m(i-1)+1,m(i)} \Vert < K_{\infty, \infty} + \varepsilon$.
Lemma \ref{estimates}d shows $Y$ $K_{\infty,\infty} + \varepsilon$-norms $X$.
Since $Y$ is independent of the blocking, and $\varepsilon$ is arbitrary,
$Y$ $c$-norms $X$.

To show the converse, we need a blocking $(m(i))$ so that the
projections $R_{m(i-1)+1,m(i)}$ are uniformly bounded. Let $m(0) = 0$,
$m(1) = 1$ and suppose $m(k)$
has been selected. Let $\varepsilon > 0$ be given. We can find unit
vectors $w_j$ which are an $\varepsilon$-net in the
sphere of $X[1,m(k)]$. We can find unit vectors $y_j$ in $Y$ so that
$|\langle w_j, y_j \rangle| > c^{-1} - \varepsilon$. Using the denseness
of $\Span\cup X^*_n$ in $Y$, select $v_j \in \Span\cup X^*_n$ so that
$\Vert v_j - y_j \Vert < \varepsilon$. Pick $m(k+1)$ so that each
$v_j \in [X^*_n]_1^{m(k+1)}$.
Let $z \in X[1,m(k)]$ have norm one, select $j$ so that
$\Vert z - w_j \Vert < \varepsilon$. We have
$$\langle z, v_j \rangle = 
	\langle w_j, y_j \rangle +
	\langle z-w_j, y_j \rangle +
	\langle z, v_j- y_j \rangle$$
$$|\langle z, v_j \rangle|  > c^{-1} - 3\varepsilon$$
And since $\langle u, v_j \rangle = 0$ if $u \in X[m(k+1)+1,\infty)$,
$\Vert z + u \Vert > (c^{-1} - 3\varepsilon)\Vert z \Vert$ and so
$\Vert R_{m(i-1)+1,m(i)} \Vert < c/(1-3\varepsilon c)$.
\end{proof}

\begin{prop}\label{anyy} If $X$ is separable and $Y$ a separating  
subspace in $X^*$, then there is a {\sf DDD} $(X_n)$ with predecomposition
space $Y$.
\end{prop}

\begin{proof}
The usual construction of a biorthogonal sequence $(x_n, x^*_n)$ (\cite{LT}
page 43), is flexible enough to produce sequences so that $X = [x_n]$
and $Y = [x^*_n]$. Then letting $X_n = [x_i]_{i=n}$, we have a {\sf DDD}
of one dimensional spaces with predecomposion space $Y$.
\end{proof}

\begin{example}\label{sepnonnorming}
There is a separating {\sf DDD} that can't be blocked to be a {\sf SBD}.
\end{example}

\begin{proof}[Construction]
Each non-quasi-reflexive space $X$ has a separable subspace $Y$ of the
dual which is separating but non-norming \cite{JL}.
\end{proof}

\begin{example}\label{onlycnorm}
There is a $c$-norming predecomposition space $Y$
that isn't one $1$-norming.
\end{example}

\begin{proof}[Construction]
Let the norm one $\phi \in X^{**}$ satisfy $0 < \dist(\phi, B_X) < 1$
where $B_X$ is  is the unit ball of $X$ as a subspace of $X^{**}$. Let
$Y = \ker \phi \subset X^*$. If $x$ has norm one and
$\Vert x - \phi \Vert = \gamma < 1$, for norm one $x^* \in Y$
$\vert x^*(x)\vert = \vert x^*(x-\phi)\vert \le \gamma < 1$, so at best
$Y$ can only $\gamma^{-1}$-norm $X$.
\end{proof}

\begin{prop}\label{renorm} If {\sf DDD} $(X_n)$ can be blocked to be
a {\sf SBD} for $X$,
then in the equivalent norm $\Vert\cdot\Vert_Y$, $(X_n)$ is still a {\sf DDD}
with the same predecompostion space $Y$ but now it is $1$-norming and
$(X_n)$ can be blocked to have asymptotic constant $1$.
\end{prop}

\begin{proof} The condition on $(X_n)$ implies $Y$ $c$-norms $X$ for some
$c < \infty$. Thus $\Vert \cdot \Vert_Y$ is an equivalent norm on $X$. 
Clearly being a {\sf DDD} or a {\sf SBD} is invariant under isomorphic
norms. Since 
$$\Vert x \Vert_Y = \lim_k \Vert q_{1,k}(x) \Vert = \lim_k\dist(x,X[k+1,\infty))$$
The quotient spaces $Z_{1,k}$ obtained for the original norm and the
new $\Vert\cdot\Vert_Y$ are isometric. It follows that the predecomposition
spaces are identical and the norms are isometric. Clearly $Y$ $1$-norms
$(X,\Vert\cdot\Vert_Y)$. The result now follows
from Theorem \ref{sbdiffnorming}.
\end{proof}

\section{Structure of {\sf SBD}  }

Basically this section shows the relationship between a {\sf SBD} $(X_n)$
and $X$ is not very strong in general. The quotient spaces $(Z_n)$
might have more importance.

\begin{prop}
\label{conv}
If $(X_n)$ is a {\sf DDD} for $X$ and $x \in X$, then there is a subsequence
$(m(i))$ and $w_i \in X[m(i-1)+1,m(i)]$ so that
$P_{m(i)} x + w_{i+1} \rightarrow x$.
\end{prop}

\begin{proof}
Let $x \in X$ be given
and suppose $(X_n)$ is {\sf DDD} and $m(k)$ has been found.
Let $M = \Vert P_{m(k)} \Vert$.
We can find $w \in \Span\cup X_n$ so that
$\Vert x - w \Vert < 2^{-k}/M $. We have
$\Vert P_{m(k)} x - P_{m(k)} w \Vert < 2^{-k}$ and
so we can let
$w_{k+1} = w - P_{m(k)} w$ and
$$\Vert x - (P_{m(k)}x + w_{k+1})\Vert
\le \Vert x - w \Vert + \Vert P_{m(k)}w - P_{m(k)}x\Vert < 2\cdot 2^{-k}.$$
Finally we select $m(k+1)$ large enough so that $w_{k+1} \in X[m(k)+1,m(k+1)]$.
\end{proof}

Example \ref{J} shows that $(\Vert  P_{m(i)} x \Vert)$ and hence
$(\Vert w_i \Vert)$ can be unbounded even when the {\sf DDD} is a {\sf SBD}.

\begin{example}
\label{J}
A {\sf SBD} of James space $J$ and a $x$ so that no subsequence of
$(P_n x)$ is bounded.
\end{example}
\begin{proof}[Construction]
James space $J$,  can be written as the set of null
sequences $(a_n)$ with finite norm given by
$\Vert (a_n) \Vert^2 = \sup \sum |a_{n(i+1)} - a_{n(i)}|^2$,
over all increasing finite integer sequences $(n(i))$. The
usual unit basis $(e_n)$ is a shrinking conditional basis which is
not boundedly complete. Let $(f_n)$ be the coefficient functionals.
Let $w_k = \sum \{e_n : 2^{k-1} < n \le 2^k\}$
and let $x = e_1+\sum_{k>0} w_k/k$. Eventually, $\Vert x \Vert = 1$.

Let $\pi$ be the permutation that reorders the integers in alternating
blocks of evens and odds in the order:
$$1,2,B_2,B_3,A_2,B_4,A_3,\dots,B_k,A_{k-1},B_{k+1},A_k\dots$$
Where
$A_k = \{n \;\mathrm{odd}: 2^{k-1} < n \le 2^k\}$ and
$B_k = \{n \;\mathrm{ even}: 2^{k-1} < n \le 2^k\}$. Let $(X_n)$ be
any blocking of $([e_{\pi(i)}])$ that is a {\sf SBD}.

Let's estimate the norm of $z_n= \sum_1^n f_{\pi(i)}(x) e_{\pi(i)}$ and
assume $\pi(n)$ is in $A_{k-1}$ or $B_{k+1}$.
This means $f_j(z_n) = 0$ for $j \in A_k$
and $f_j(z_n) = 1/k$ for $j \in B_k$. Thus
$\Vert z_n \Vert \ge \sqrt{2^k/k^2}$. It follows that every subsequence
of $(P_n x)$ is unbounded no matter how $([e_{\pi(i)}])$ was blocked.
\end{proof}

\begin{lemma}
If $(X_n)$ is a {\sf DDD} and $x\in X$, then $\lim q_n(x) = 0$
\end{lemma}

\begin{proof}
Let $x \in X$ and $\varepsilon > 0$.  Since span $(X_n)$ is dense in $X$
we can find $N$ and $y \in X[1,N]$ so that $\Vert x-y \Vert < \varepsilon$.
Since $x-(x-y) = y$ is in $X[1,N]$,
$q_n(x) = q_n(x-y)$ and $\Vert q_n(x-y)\Vert  < \varepsilon$
for $n>N$.
\end{proof}

\begin{example}
\label{J2}
A {\sf SBD} $(X_n)$ and $x \in X$ so that $(p_n(x))$ has no bounded
subsequence. In particular, the isomorphisms $q_n |_{X_n} : X_n \rightarrow
Z_n$ have inverses with norms that blow up. Thus $X_n$ and $Z_n$ can
be far apart.
\end{example}
\begin{proof}[Construction]
This is a continuation of Example \ref{J}. For $k > 3$, let $X_k =
[e_i: i \in B_k \cup A_{k-1}]$. The same $x$, has $\Vert p_k(x) \Vert
\ge \sqrt{2^k/k^2}$, while $q_k(x) = q_k(p_k(x)) \rightarrow 0$ by the
lemma.
\end{proof}

\section {$Y^*$ and Boundedly Complete {\sf SBD}}

\begin{define}
\label{bcsbddn}
If $x_n \in X_n$ and $(\Vert \sum_1^n, x_i \Vert)$ bounded implies that
$\sum x_i$ converges in norm, then we call $(X_n)$ is boundedly complete.
A boundedly complete {\sf SBD} is one where all the skipped decompositions
(X[n(i),m(i)]) in equation (3) are boundedly complete. (Note that the whole
sequence $(X_n)$ is not required to be boundedly complete.)
\end{define}

The subspace map $\phi: Y \rightarrow X^*$ yields by duality
a quotient map $\phi^* : X^{**} \rightarrow Y^*$. The restriction
of $\phi^*$ to the image of $X$ in $X^{**}$ is the duality given by
$\langle y, \phi^* x \rangle = \langle x, y \rangle$ for $y \in Y$.
The adjoint of quotient maps $q_{m,k}: X \rightarrow Z_{m,k}$ given by (7),
factors through $Y$ as the subspace injections
$Z^*_{m,k} \subset Y \subset X^*$
so the double adjoint quotient map
$q^{**}_{m,k}: X^{**} \rightarrow Z^{**}_{m,k} = Z_{m,k}$ factors 
through $Y^*$. Similarly we can identify $X_n$ in $Y^*$ (as a set
with perhaps a different norm) as the range of the projection $p_n^{**}
=\phi^*(X_n)$. Note that $\phi^*$ is an isomorphism exactly when $Y$
is norming. Theorem \ref{sbdiffnorming} implies $\phi^*$ is an
isomorphism exactly when the {\sf DDD} $(X_n)$  can be blocked to
be a {\sf SBD}. Note further that $\phi^*$ is an isometry if $X$
has {\sf PCP} or if $(X_n)$ is shrinking.

On the unit ball, $B_{Y^*}$, of $Y^*$ there are several differently
defined weak topologies which are the same. There is
$\sigma(X^{**},X^{*})$-topology on $B_{X^{**}}$ which the quotient
map $\phi^*$ induces on $B_{Y^*}$ which must be the $\sigma(Y^*,Y)$-
topology by compactness and weak-star continuity. Also this is
the same as the $\sigma(Y^*, \cup Y_n)$, also because of compactness.
In this topology, $B_{Y^*}$ is a compact metric space.

For $G$ a finite subset of $\cup Y_n$
and $\varepsilon > 0$, let $$V(x, G, \varepsilon) = \{y^* \in Y^{*}:
|<y^*-x,g>| < \varepsilon, \forall g \in G\}.$$
it follows that the collection of $V(x,G,\varepsilon)$ is a basis
for the $\sigma(Y^*,Y)$-topology on $Y^*$. We use the term
{\em elementary } to describe $\sigma(Y^*,Y)$-open sets of
the form $V(x,G,\varepsilon)$ with $G \subset \cup Y_n$.

\begin{thm}
\label{allbc}
If $X$ has {\sf PCP} then each norming {\sf DDD}
can be blocked to
be a boundedly complete {\sf SBD}.
\end{thm}

\begin{proof}
Let $(X_n)$ be a {\sf DDD} with a $c$-norming predecomposition space
for $X$. Renorming via Proposition \ref{renorm} if necessary, we can assume
the predecomposition space $Y$ is $1$-norming. By blocking if necessary
 we can assume that
$(X_n)$ is already a {\sf SBD}. We pick $(m(k))$ by
induction so that $m(0)=0$ and $(X[m(k-1)+1,m(k)])$ is a boundedly complete
{\sf SBD}. Our mode of proof is to make slight modifications to
the existence proof in \cite{GM1}.

The {\sf PCP} property (Lemma II.3 \cite{GM1}) implies that
$Y^*\backslash X = \cup K_n$ for some increasing sequence 
$(K_n)$ of $\sigma(Y^*,Y)$-compact sets. While the cited lemma
constructs a particular $Y$, any $1$-norming subspace of $X^*$
will satisfy the proof.

Suppose $m(k)$ has been selected. For each $x$ in the unit ball
of $X[1,m(k)]$ there is an elementary $\sigma(Y^*,\cup Y_n)$-open
$V$ with $V \cap K_{k+1} = \emptyset$. By compactness
we can assume only a finite number of the open sets
$V_i = V(x_i,G_i,\varepsilon_i)$ are needed.
Let $m(k+1)$ be large enough so that $Y[1,m(k+1)]$ contains all
the vectors in $G_i$ and also $1+1/k$-norms $X[1,m(k)]$. 

The claim is that the {\sf SBD} $(X[m(k-1)+1,m(k)])_k$ is also a
boundedly complete {\sf SBD}. Suppose $x_k \in X[m(k-1),m(k)]$
with $x_k = 0$ infinity often. Let $n(i)$ be the sequence defined
so that $x_{n(i)} \not= 0$
but $x_{n(i)+1} = 0$ and $s_i =\sum_1^{n(i)} x_k$ be so that
$\Vert s_i\Vert \le 1$ uniformly. Let $s \in Y^*$ be the
$\sigma(Y^*,Y)$-limit of $(s_i)$.

If $s \in Y^* \backslash X$ then $s \in K_j$ for $j\ge N$ for some $N$. Let
$k=n(i) > N$ so that $x_{k+1} = 0$. Then there is $V=V(w,G,\varepsilon)$
so that $s_i \in V$ and $V \cap K_{k+1} = \emptyset$.
Since $G \subset Y[1,m(k+1)]$, 
$<x_i, g> = 0$ for $i > k+1$. Since $x_{k+1} = 0$ This means $s_j \in V$
for all $j \ge i$, so $s \in V \cap K_{k+1}$ a contradiction.

Therefore $s \in X$, and we need to show $s_i$ converges in norm to $s$.
Let $\varepsilon  > 0$ be given. Find $k = n(j)$ and $w \in X[1,m(k)]$ so that
$\Vert s - w \Vert < \varepsilon$ and $x_{k+1} = 0$. For $i \ge j$
Find $y \in Y[1,m(n(j)+1)]$
with $\Vert y \Vert = 1$ and $y$ almost norms $s_i - w$. In particular
$\langle y, s_i - w \rangle \le \Vert s_i - w \Vert(1+1/i)^{-1}$. 
Since $\langle y, x_n \rangle = 0$ for $n > n(i)$, we have
$\langle y, s-w \rangle = \langle y, s_i - w \rangle$. Thus
$\Vert s_i - w \Vert \le \Vert s - w \Vert (1 + 1/i)$ and
$\Vert s_i - s \Vert \le 3\varepsilon$. Thus $(s_i)$ norm converges to $s$.
\end{proof}

\section{{\sf JT}-decompositions}

\begin{define}
\label{JTdn}
The boundedly complete {\sf SBD} $(X_n)$ is called a {\em 
$c$-norming-{\sf JT}-decomposition} for $X$ provided
\begin{enumerate}
\item[(A)]For each $x^{**}\in X^{**}$, with $\liminf\Vert 
q^{**}_n(x^{**})\Vert=0$, there is $x\in X$ with $\Vert x\Vert\leq 
c\Vert x^{**}\Vert$ and for all $n$, $p^{**}_n(x^{**})=p_n(x)$. 
\end{enumerate}
\end{define}

Ghoussoub, Maurey and Schachermayer \cite{GMS} have defined a 
{\sf JT}-decomposition to be equivalent to what we have called a 
$1$-norming-{\sf JT}-decomposition.  Furthermore they show each separable 
Banach space with {\sf PCP} has a $1$-norming-{\sf JT}-decomposition.  The 
following proposition shows every boundedly complete {\sf SBD} is a 
$c$-norming-{\sf JT}-decomposition.

\begin{prop}
\label{bcdJT}
{If $(X_n)$ is a boundedly complete {\sf SBD} for $X$, 
then $(X_n)$ is a $K_{\infty}$-norming-{\sf JT}-decomposition where 
$K_{\infty}$ is the asymptotic-{\sf SBD}-constant of $(X_n)$.}
\end{prop}

\begin{proof}
First note that if $x\in X\backslash\{0\}$, then
Lemma \ref{estimates}c implies 
for some $n$, $p_n(x)\not=0$.  It follows that if such an $x\in X$ 
as in (A) exists, then it must be unique.

Let $x^{**}\in X^{**}$ and let $\varepsilon>0$ be given.  Find $m(1)$ so 
that $\Vert q^{**}_{m(1)}(x^{**})\Vert<\varepsilon/2^1$ and let 
$w_1\in\Span\cup X_n$ so that 
$q_{m(1)}(w_1)=q^{**}_{m(1)}(x^{**})$ and $\Vert 
w_1\Vert<\varepsilon/2^1$.  Suppose $(w_i)^k_{i=1}$ and $(m(i))^k_{i=1}$ 
have been chosen so that

\begin{enumerate}
\item[(i)]$q_{m(i)}(w_i)=q^{**}_{m(i)}(x^{**})$;

\item[(ii)]$\Vert q_{m(i)}(w_i)\Vert<\varepsilon/2^i$;

\item[(iii)]$q_{m(j)}(w_i)=0$ for $j<i\leq k$;

\item[(iv)]$w_{i-1}\in X[1,m(i)-1]$ for 
$i\leq k$; and

\item[(v)]$w_k\in\Span\cup X_n$.
\end{enumerate}

Let $C=\Vert I-p_{m(1)}\Vert\;\Vert I-p_{m(2)}\Vert\cdots\Vert I-
p_{m(k)}\Vert$ and pick $m(k+1)$ large so that
$w_k\in X[1,m(k+1)-1]$ and $\Vert 
q^{**}_{m(k+1)}(x^{**})\Vert<\varepsilon/C2^{k+1}$.  Find $z\in\Span\cup 
X_n$ so that $\Vert z\Vert<\varepsilon/C2^{k+1}$ and 
$q_{m(k+1)}(z)=q^{**}_{m(k+1)}(x^{**})$.  Let $w_{k+1}=(I-
p_{m(1)})(I-p_{m(2)})\cdots(I-p_{m(k)})z$.  Clearly $\Vert 
w_{k+1}\Vert<\varepsilon/2^{k+1}$ and $q_{k(i)}(w_{k+1})=0$ for $i\leq k$.  
Therefore, writing $w = \sum w_k$,

\begin{enumerate}
\item[(B)]For each $\varepsilon>0$ there is $w\in X$ and $(m(i))$ 
so that $\Vert w\Vert<\varepsilon$ and $q_{m(i)}(x^{**}-w)=0$ for all $i$.
\end{enumerate}

Now for each $i$, $x^{**}-
w\in[X_n]^{\perp\perp}_{n\not=m(i)}=X^{*\perp}_{m(i)}$, hence 
$$\Vert \sum\nolimits^{m(i)-1}_{j=1}p^{**}_j(x^{**}-w)\Vert=\Vert 
R_{m(i)}(x^{**}-w)\Vert \leq\Vert R_{m(i)}\Vert\;\Vert x^{**}-
w\Vert\leq K\Vert x^{**}-w\Vert.$$  Since $(X[m(i-1)+1,m(i)-1])_i$ 
is boundedly complete 
there is $z\in X$ with
$$\Vert 
z\Vert\leq\lim\sup\Vert R_{m(i)}\Vert\;\Vert x^{**}-w\Vert\leq 
K_{\infty}\Vert x^{**}-w\Vert$$ and $p_n(z)=p^{**}_n(x^{**}-w)$ for 
all $n$.  Finally let $x=z+w$.  We have $p_n(x)=p^{**}_n(x^{**})$ 
and $$\Vert x\Vert=\Vert z+w\Vert<\Vert z\Vert+\varepsilon\leq 
K_{\infty}\Vert x^{**}-w\Vert+\varepsilon\leq K_{\infty}\Vert 
x^{**}\Vert+(K_{\infty}+1)\varepsilon.$$  The choice of $x$ is independent 
of $\varepsilon$ and so $\Vert x\Vert\leq K_{\infty}\Vert 
x^{**}\Vert$.
\end{proof}

The following theorem summarizes our results for Banach space with {\sf PCP}.

\begin{thm}\label{climax}

The following are equivalent for a Banach space $X$:
\begin{enumerate}
\item[(a)]$X$ is separable and has {\sf PCP}.

\item[(b)]Each norming {\sf DDD} $(X_n)$ of $X$ can be blocked to be a 
boundedly complete {\sf SBD}.

\item[(c)]Each norming {\sf DDD} $(X_n)$ of $X$ can be blocked to be a 
$1$-norming {\sf JT}-decomposition in an equivalent norm.
\end{enumerate}
\end{thm}

\begin{proof}
(c)$\Rightarrow$(b) is formal. (b)$\Rightarrow$(a) is in
Bourgain and Rosenthal \cite{BR}. (a)$\Rightarrow$(b) is Theorem \ref{allbc}.
(b)$\Rightarrow$(c) uses Proprosition \ref{renorm}
to equivalently renorm $X$ so that the predecomposition
space is $1$-norming, Theorem \ref{sbdiffnorming} to block the {\sf DDD}
to a {\sf SBD} with asymptotic constant $1$ and finally Theorem \ref{bcdJT}
to show this blocking is a $1$-norming {\sf JT}-decompostion.
\end{proof}

{\sc
Department of Mathematics, Florida State University, Tallahassee, FL 32306-4510
}

{\em E-mail address:}\;{\bf bellenot@math.fsu.edu}


\begin{thebibliography}{WWW}

\bibitem{B1} S.F. Bellenot {\em Somewhat quasireflexive Banach spaces}
Ark. Mat. 22(1984), 175-183.

\bibitem{B2} S.F. Bellenot {\em More quasi-reflexive subspaces}
Proc. Amer. Math. Soc. 101(1987), 693-696.

\bibitem{BR} J. Bourgain and H.P. Rosenthal, {\em Geometrical 
implications of certain finite-dimensional decompositions}, Bull. 
Soc. Math. Belg. 32(1980), 57-82.

\bibitem{DL} W.J. Davis and J. Lindenstrauss, {\em On total nonnorming
subspaces}, Proc. Amer. Math. Soc 31 (1972), 109-111.

\bibitem{EW} G.A. Edgar and R.F. Wheeler, {\em Topological 
properties of Banach spaces}, Pacific J. Math. 115(1984), 317-350.

\bibitem{F} C. Finet, {\em Subspaces of Asplund Banach spaces with 
the point of continuity property}, Israel J. Math. 60(1987), 
191-198.

\bibitem{GM1} N. Ghoussoub and B. Maurey, {\em $G_\delta$-embeddings in 
Hilbert space}, J. Funct. Anal 61 (1985) 72-97.

\bibitem{GM2} N. Ghoussoub and B. Maurey, {\em $H_\delta$-embeddings in 
Hilbert space}, J. Funct. Anal 78 (1988) 271-305.

\bibitem{GMS} N. Ghoussoub, B. Maurey and W. Schachermayer, {\em 
Geometrical implications of certain infinite dimensional 
decompositions}, Trans. Amer. Math. Soc. 317(1990), 541-584.

\bibitem {JL} W.B. Johnson and J. Lindenstrauss, {\em Basic concepts in
the geometry of Banach space}, in \cite{JLE} (2001), 1-84.

\bibitem {JLE} W.B. Johnson and J. Lindenstrauss (Editors),
{\em Handbook of the geometry of Banach spaces, Volume 1}
Elsevier, Amsterdam, 2001.

\bibitem {JR} W.B. Johnson and H.P. Rosenthal, {\em On $w^*$-basic sequences and
their applications to the study of Banach spaces}, Studia Math 43 (1972),
77-92.

\bibitem{LT} J. Lindenstraus and L. Tzafriri, {\em Classical Banach
spaces I: Sequence spaces} Springer-Verlag, Berlin and New York, 1977.

\bibitem{R} H.P. Rosenthal, {\em Weak$^*$-polish Banach spaces}, J. 
Funct. Anal. 76(1988), 267-316.
\end{thebibliography}
\end{document}